\documentstyle[11pt]{article}

\makeatletter

\@addtoreset{equation}{section}
\makeatother
\def\<{\langle}
\def\>{\rangle}

\newtheorem{lem}{Lemma}[section]
\newtheorem{theo}{Theorem}[section]
\newtheorem{rem}{Remark}[section]

\makeatletter
   
   \@addtoreset{equation}{section}
\makeatother

\begin{document}
\title{\bf Asymptotic Profiles for Wave Equations\\ with Strong Damping}
\author{Ryo IKEHATA\thanks{Corresponding author: ikehatar@hiroshima-u.ac.jp} \\ {\small Department of Mathematics, Graduate School of Education, Hiroshima University} \\ {\small Higashi-Hiroshima 739-8524, Japan}}
\maketitle
\begin{abstract}
We consider the Cauchy problem in ${\bf R}^{n}$ for strongly damped wave equations. We derive asymptotic profiles of these solutions with weighted $L^{1,1}({\bf R}^{n})$ data by using a method introduced in \cite{INa}. 
\end{abstract}

\section{Introduction}
\footnote[0]{Keywords and Phrases: Wave equation; Strong damping; Fourier analysis; Asymptotic profiles; Low frequency; Weighted $L^{1}$-initial data.}
\footnote[0]{2010 Mathematics Subject Classification. Primary 35B40, 35L15; Secondary 35L05, 35C20.}
We are concerned with the Cauchy problem for strongly damped wave equations in ${\bf R}^{n}$ ($n \geq 1$):
\begin{equation}
u_{tt}(t,x) - \Delta u(t,x)  -\Delta u_{t}(t,x) = 0,\ \ \ (t,x)\in (0,\infty)\times {\bf R}^{n},
\end{equation}
\begin{equation}
u(0,x)= u_{0}(x), \quad u_{t}(0,x) = u_{1}(x), \quad x\in {\bf R}^{n},
\end{equation}
where the initial data $u_{0}$ and $u_{1}$ are taken from the energy space:
\[[u_{0},u_{1}] \in H^{1}({\bf R}^{n})\times L^{2}({\bf R}^{n}).\]
\noindent
It is known (see \cite{ITY}) that the problem (1.1)-(1.2) admits a unique weak solution\\
 $u \in C([0,+\infty);H^{1}({\bf R}^{n})) \cap C^{1}([0,+\infty);L^{2}({\bf R}^{n}))$.\\

The purpose of this paper is to investigate the asymptotic profiles of the solution $u(t,x)$ to problem (1.1)-(1.2). As for the related results Ponce \cite{p} and Shibata \cite{shibata} have already studied the $L^{p}-L^{q}$ decay estimates of solutions to problem (1.1)-(1.2). The exterior domain case for the equation (1.1) was also handled in Ikehata \cite{Ik-1}, in which the two dimensional case is likely to be sharp about the rate of decay of the corresponding total energy and the $L^{2}$-norm of solutions. Recently Ikehata-Natsume \cite{INatsume} derived the decay estimates of the total energy and $L^{2}$-norm of solutions to (1.1)-(1.2) with a more general structural damping based on the energy method in the Fourier space due to \cite{UKS}. Before \cite{INatsume}, Ikehata-Todorova-Yordanov \cite{ITY} succeeded to find the  asymptotic profile in the abstract framework, in fact, they studied the following OD equations in Hilbert space $H$:
\begin{equation}
u_{tt}(t) + Au(t) + Au_{t} = 0,\quad u(0) = u_{0}, \quad u_{t}(0) = u_{1},
\end{equation}
where $A$ is a nonnegative self-adjoint operator in $H$. They employed the abstract energy method in the Fourier space combined with the spectral analysis to find the asymptotic profile such as 
\[u(t) \sim e^{-tA/2}(\cos(tA^{1/2})u_{0} + A^{-1/2}\sin(tA^{1/2})u_{1}), \quad (t \to +\infty).\]
These ideas are inspired from \cite{CH} and \cite{UKS}. Therefore, in this sense we have already caught the asymptotic profiles from the work \cite{ITY}, however, it seems to be important to search another root to find the profiles of the solution to (1.1)-(1.2) by the concrete way because one can sometimes find a possibility of several new applications of the method introduced in this paper. Our new point of view is in dealing with the problem (1.1)-(1.2) in a framework of  the weighted $L^{1}$-data. By imposing some weights on the initial data in $L^{1}$ sense we can get the meaningful "equality" represented by (2.11) below, which includes explicitly the leading plus error terms. This method is basically independent from the shape of equation itself. The origin of this idea comes from \cite{Ik-2, Ik-3}, which studied the decay property of solutions to the damped wave equations: 
\begin{equation}
u_{tt}(t,x) - \Delta u(t,x) + u_{t}(t,x) = 0.
\end{equation}
The research in the framework of the weighted $L^{1}$-data was also developed more precisely to (1.4) in the recent results due to \cite{KU} and \cite{INa}. Especially \cite{KU} dealt with the nonlinear problems of (1.4). The decay property and the asymptotic profiles to the equation (1.4) are well-studied in \cite{HM}, \cite{INa}, \cite{IN}, \cite{K}, \cite{Ma}, \cite{Na}, \cite{N-2} and \cite{RTY}. As compared with the equation (1.4), there seems to be few results about the asymptotic profiles of solutions to the equation (1.1), so it is good chance to present a way to investigate the asymptotic behavior of solutions to (1.1)-(1.2). In this connection, quite recently vigorous works about the global existence of solutions and/or a new method to derive sharper decay estimates of the total energy to the Cauchy problem of the equation
\[u_{tt} - \Delta u + (-\Delta)^{\sigma}u_{t} = \mu f(u),\quad \sigma \in [0,1]\]
are successively announced by D'Abbicco-Reissig \cite{DR} ($\mu > 0$) and Char\~ao-daLuz-Ikehata \cite{RCR} ($\mu = 0$), respectively.\\  

Our main target is to obtain the asymptotic profile of the solution $u(t,x)$ to problem (1.1)-(1.2) as $\to +\infty$ as follows.
\begin{theo}\,Let $n \geq 1$. If $[u_{0},u_{1}] \in (H^{1}({\bf R}^{n}) \cap L^{1,1}({\bf R}^{n})) \times (L^{2}({\bf R}^{n}) \cap L^{1,1}({\bf R}^{n}))$, then the solution $u(t,x)$ to problem {\rm (1.1)}-{\rm (1.2)} satisfies
\[\int_{{\bf R}^{n}}\vert{\cal F}(u(t,\cdot))(\xi)- \{P_{1}e^{-t\vert\xi\vert^{2}/2}\frac{\sin(t\vert\xi\vert)}{\vert\xi\vert}+P_{0}e^{-t\vert\xi\vert^{2}/2}\cos(t\vert\xi\vert)\}\vert^{2}d\xi\]
\[\leq C(\Vert u_{1}\Vert_{1}^{2} + \Vert u_{0}\Vert_{1}^{2})t^{-\frac{n}{2}} + C\Vert u_{1}\Vert_{1,1}^{2}t^{-\frac{n}{2}} + C\Vert u_{0}\Vert_{1,1}^{2}t^{-\frac{n}{2}-1} + Ce^{-\alpha t}(\Vert u_{1}\Vert^{2} + \Vert u_{0}\Vert^{2})\]
with some constants $\alpha > 0$ and $C > 0$, where
\[P_{0} := \int_{{\bf R}^{n}} u_{0}(x)dx,\quad P_{1} := \int_{{\bf R}^{n}} u_{1}(x)dx.\]
\end{theo}
\begin{rem} {\rm It follows from \cite[(i) of Theorem 1.4 with $\theta = 1$]{INatsume} that even if we assume $[u_{0},u_{1}] \in L^{1,1}({\bf R}^{n})\times L^{1,1}({\bf R}^{n})$,  we find that
\[\Vert u(t,\cdot)\Vert^{2} \leq C(1+t)^{-\frac{n}{2}}\Vert u_{1}\Vert_{1,1}^{2} +  C(1+t)^{-\frac{n+2}{2}}\Vert u_{0}\Vert_{1,1}^{2} \]
\[+ Ce^{-\eta t}(\Vert u_{0}\Vert^{2}+\Vert u_{1}\Vert^{2} ) + C(1+t)^{-\frac{n-2}{2}}\vert P_{1}\vert^{2} + C(1+t)^{-\frac{n}{2}}\vert P_{0}\vert^{2} ,\]
provided that $n \geq 3$.  This implies that in the case of $P_{1} \ne 0$, we have at most $\Vert u(t,\cdot)\Vert^{2} = O(t^{-\frac{n-2}{2}})$ as $t \to +\infty$. On the other hand, we can observe by a simple computation that
\[\int_{{\bf R}^{n}}e^{-t\vert\xi\vert^{2}}\vert\frac{\sin(t\vert\xi\vert)}{\vert\xi\vert}\vert^{2}d\xi = O(t^{-\frac{n-2}{2}}), \quad (t \to +\infty),\]
if $n \geq 3$. Furthermore, in the case when $n \geq 1$ we see that
\[\int_{{\bf R}^{n}}e^{-t\vert\xi\vert^{2}}\vert\cos(t\vert\xi\vert)\vert^{2}d\xi = O(t^{-\frac{n}{2}}), \quad (t \to +\infty).\]
This observation implies that the asymptotic profile as $t \to +\infty$ of the solution $u(t,x)$ to problem (1.1)-(1.2) becomes
\begin{equation}
{\cal F}^{-1}(e^{-t\vert\xi\vert^{2}/2}\{\frac{\sin(t\vert\xi\vert)}{\vert\xi\vert}P_{1} + \cos(t\vert\xi\vert)P_{0}\})(x)
\end{equation}
at least formally in the case when $P_{1} \ne 0$ and $n \geq 3$.}
\end{rem}
\begin{rem} {\rm We give a remark to the case of  $n = 2$. It follows from \cite[(3) of Theorem 2.1]{shibata} that if $[u_{0},u_{1}] \in L^{1}({\bf R}^{2})\times L^{1}({\bf R}^{2})$,  then we find that
\[\Vert u(t,\cdot)\Vert^{2} = O(\{\log(1+ t)\}^{2})\quad (t \to +\infty),\]
provided that $\Vert u_{1}\Vert_{1} \ne 0$, $n = 2$.  On the other hand, if we set $M := \displaystyle{\sup_{\theta \ne 0}}\{\vert\sin \theta\vert/\vert\theta\vert\}$,
then it follows that
\[\int_{{\bf R}^{2}}e^{-t\vert\xi\vert^{2}}\vert\frac{\sin(t\vert\xi\vert)}{\vert\xi\vert}\vert^{2}d\xi = t^{2}\int_{{\bf R}^{2}}e^{-t\vert\xi\vert^{2}}\vert\frac{\sin(t\vert\xi\vert)}{t\vert\xi\vert}\vert^{2}d\xi \]
\[\leq M^{2}t^{2}\int_{{\bf R}^{2}}e^{-t\vert\xi\vert^{2}}d\xi = O(t), \quad (t \to +\infty).\]
From this observation we can find that the asymptotic profile of the solution $u(t,x)$ to problem (1.1)-(1.2) with $n = 2$ becomes the same as Remark 1.1. On the other hand, as for the one dimensional case ($n = 1$) we do not have any previous knowledges about the $L^{2}$-decay rate of the solution $u(t,\cdot)$ to (1.1)-(1.2), we can say nothing at present, however, the situation will be the same as the case when $n \geq 2$.}
\end{rem}

Final part of this section is devoted to represent explicitly the formula (1.5) based on the well-known fact called as the Kirchhoff formulas for solutions of the free wave equation (see Evans \cite{E} and Mizohata \cite{M} and Shibata \cite{shibata}) :
\begin{equation}
w_{tt}(t,x) - \Delta w(t,x) = 0,\ \ \ (t,x)\in (0,\infty)\times {\bf R}^{n},
\end{equation}
\begin{equation}
w(0,x)= 0, \quad w_{t}(0,x) = \delta, \quad x\in {\bf R}^{n},
\end{equation}
where $\delta(x)$ is the usual Dirac measure. It is well-known that the Fourier images of $w$ and $w_{t}$ are given by
\[\hat{w}(t,\xi) = \frac{\sin(t\vert\xi\vert)}{\vert\xi\vert}, \quad\hat{w}_{t}(t,\xi) = \cos(t\vert\xi\vert),\]
respectively.
\begin{lem}\,Let $n \geq 2$. Then there are constants $a_{\alpha}$ and $b_{\alpha}$ such that if $w(t,x)$ is the fundamental solution of the free wave equation {\rm (1.6)}-{\rm (1.7)}, and if $h \in C^{\infty}({\bf R}^{n})$, then
in the case when $n$ is odd:
\begin{equation}
(w\ast h)(t,x) = \sum_{0 \leq \vert\alpha\vert\leq (n-3)/2}a_{\alpha}t^{1+\vert\alpha\vert}\int_{\vert z\vert = 1}z^{\alpha}(\partial_{x}^{\alpha}h)(x+t z)dS_{z},
\end{equation}
\begin{equation}
(w_{t}\ast h)(t,x) = \sum_{0 \leq \vert\alpha\vert\leq (n-1)/2}b_{\alpha}t^{\vert\alpha\vert}\int_{\vert z\vert = 1}z^{\alpha}(\partial_{x}^{\alpha}h)(x+t z)dS_{z},
\end{equation}
and in the case where $n$ is even:
\begin{equation}
(w\ast h)(t,x) = \sum_{0 \leq \vert\alpha\vert\leq (n-2)/2}a_{\alpha}t^{1+\vert\alpha\vert}\int_{\vert z\vert \leq 1}\frac{z^{\alpha}(\partial_{x}^{\alpha}h)(x+t z)}{\sqrt{1-\vert z\vert^{2}}}dz,
\end{equation}
\begin{equation}
(w_{t}\ast h)(t,x) = \sum_{0 \leq \vert\alpha\vert\leq n/2}b_{\alpha}t^{\vert\alpha\vert}\int_{\vert z\vert \leq 1}\frac{z^{\alpha}(\partial_{x}^{\alpha}h)(x+t z)}{\sqrt{1-\vert z\vert^{2}}}dz,
\end{equation}
where $dS_{z}$ denotes surface measure on the unit sphere in ${\bf R}^{n}$ and $\alpha := (\alpha_{1},\alpha_{2},\cdots,\alpha_{n})$ denotes a multi-indices.
\end{lem} 
Let
\[G(t,x) := C_{n}t^{-\frac{n}{2}}e^{-\frac{\vert x\vert^{2}}{2t}} = {\cal F}^{-1}(e^{-t\vert\xi\vert^{2}/2}),\]
with some constant $C_{n} > 0$ depending on each $n$. Then by Lemma 1.1, in the case when $n$ is odd one has
\[{\cal F}^{-1}(e^{-t\vert\xi\vert^{2}/2}\frac{\sin(t\vert\xi\vert)}{\vert\xi\vert}P_{1})(x) = P_{1}\cdot(w(t,\cdot)\ast G(t,\cdot))(x)\]
\[= \sum_{0 \leq \vert\alpha\vert\leq (n-3)/2}a_{\alpha}\,t^{1+\vert\alpha\vert}\int_{\vert z\vert = 1}z^{\alpha}(\partial_{x}^{\alpha}G)(t,x + tz)dS_{z}\cdot P_{1},\]
and 
\[{\cal F}^{-1}(e^{-t\vert\xi\vert^{2}/2}\cos(t\vert\xi\vert)P_{0})(x) = P_{0}\cdot(w_{t}(t,\cdot)\ast G(t,\cdot))(x)\]
\[= \sum_{0 \leq \vert\alpha\vert\leq (n-1)/2}b_{\alpha}\,t^{\vert\alpha\vert}\int_{\vert z\vert = 1}z^{\alpha}(\partial_{x}^{\alpha}G)(t,x + tz)dS_{z}\cdot P_{0}.\]
As a result the asymptotic profile of the solution $u(t,x)$ as $t \to +\infty$ to problem (1.1)-(1.2) can be written exactly by the convolution of the fundamental  solutions of diffusion and (free) wave equations. We write down the obtained formulas symbolically:
\[u(t,x) \sim Ct^{-1/2}(\int_{{\bf R}^{3}}u_{1}(y)dy)\int_{\vert z\vert = 1}e^{-\frac{\vert x + tz\vert^{2}}{2t}}dS_{z},\quad (t \to +\infty) \hspace{0.2cm}\textstyle{in}\hspace{0.2cm}L^{2}({\bf R}^{3}),\]
in the case when $P_{0} = 0$ (for simplicity) and $n = 3$, and
\[u(t,x) \sim C(\int_{{\bf R}^{2}}u_{1}(y)dy)\int_{\vert z\vert \leq 1}\frac{e^{-\frac{\vert x + tz\vert^{2}}{2t}}}{\sqrt{1-\vert z\vert^{2}}}dz,\quad (t \to +\infty) \hspace{0.2cm}\textstyle{in}\hspace{0.2cm}L^{2}({\bf R}^{2}),\]
in the case when $P_{0} = 0$ (for simplicity) and $n = 2$. By using (1.9) and (1.11) one can also get the explicit formulas in the case of nontrivial $P_{0} \ne 0$. \\
\noindent
In any case we can find that the solution $u(t,x)$ is time-asymptotic to the flow $w\ast (GP_{1}) + w_{t}\ast (GP_{0})$ with total mass $P_{j}$ (j = 0,1), which is closely related with the diffusion waves in the field of the Navier-Stokes equations of compressible flow. Here, $G(t,\cdot)P_{j}$ is the solution of some parabolic equations with initial datum $P_{j}\delta$ ($j = 0,1$). The diffusion waves in the Navier-Stokes equations are well-studied in Hoff-Zumbrun \cite{HZ} and the references therein, and in particular, the weighted $L^{1}$-initial data (i.e., $L^{1,1}$ initial data) have been already assumed to capture the diffusion wave property in \cite[(1.9) of Theorem, Theorem 6.8]{HZ} from a little different viewpoint. As far as we know there seems no any previous works attacking to the equation (1.1) directly to investigate the asymptotic profiles except for \cite{ITY}.\\

{\footnotesize {\bf Notation.} Throughout this paper, $\| \cdot\|_q$ stands for the usual $L^q({\bf R}^{n})$-norm. For simplicity of notations, in particular, we use $\| \cdot\|$ instead of $\| \cdot\|_2$. 
Furthermore, we set
\[f \in L^{1,\gamma}({\bf R}^{n}) \Leftrightarrow f \in L^{1}({\bf R}^{n}), \Vert f\Vert_{1,\gamma} := \int_{{\bf R}^{n}}(1+\vert x\vert)^{\gamma}\vert f(x)\vert dx < +\infty, \quad \gamma \geq 0.\]\\
On the other hand, we denote the Fourier transform $\hat{\phi}(\xi)$ of the function $\phi(x)$ by
\begin{equation}
{\cal F}(\phi)(\xi) := \hat{\phi}(\xi) := \frac{1}{(2\pi)^{n/2}}\int_{{\bf R}^{n}}e^{-ix\cdot\xi}\phi(x)dx,
\end{equation}
where $i := \sqrt{-1}$, and $x\cdot\xi = \displaystyle{\sum_{i=1}^{n}}x_{i}\xi_{i}$ for $x = (x_{1},\cdots,x_{n})$ and $\xi = (\xi_{1},\cdots,\xi_{n})$, and the inverse Fourier transform of ${\cal F}$ is denoted by ${\cal F}^{-1}$. When we estimate several functions by applying the Fourier transform sometimes we can also use the following definition in place of (1.12)
\[{\cal F}(\phi)(\xi) := \int_{{\bf R}^{n}}e^{-ix\cdot\xi}\phi(x)dx\]
without loss of generality. We also use the notation
\[v_{t}=\frac{\partial u}{\partial t}, \quad v_{tt}=\frac{\partial^{2} v}{\partial t^{2}}, \quad \Delta = \sum^n_{i=1}\frac{\partial^2}{\partial x_i^2},\ \ x=(x_1,\cdots,x_n), \quad (f\ast g)(x) := \int_{{\bf R}^{n}}f(x-y)g(y)dy.\]}


\section{Proof of Theorem 1.1.}

Let us prove Theorem 1.1 based on an idea due to \cite{INa} which has its origin in \cite{Ik-3}. The essential part of the proof corresponds to the low frequency estimates of the solution. For the moment, we shall assume that the initial data $[u_{0},u_{1}]$ are sufficiently smooth.\\

Let $\hat{u}(t,\xi) := {\cal F}(u(t,\cdot))(\xi)$. Then we first prove the following lemma.
\begin{lem}  Let $n \geq 1$. Then, it is true that there exists a constant $C>0$ such that for $t > 0$ 
\[\int_{\vert\xi\vert\leq \delta_{0}}\vert \hat{u}(t,\xi) - \{P_{1}e^{-t\vert\xi\vert^{2}/2}\frac{\sin(t\vert\xi\vert)}{\vert\xi\vert} + e^{-t\vert\xi\vert^{2}/2}P_{0}\cos(t\vert\xi\vert)\}\vert^{2}\,d\xi\]
\[\leq Ct^{-\frac{n}{2}-1}\Vert u_{0}\Vert_{1,1}^{2}+Ct^{-\frac{n}{2}}\Vert u_{1}\Vert_{1,1}^{2}\]
for small positive $\delta_{0} \ll 1$.
\end{lem}
To begin with, we apply the Fourier transform with respect to the space variable $x$ of both sides of (1.1)-(1.2). Then in the Fourier space ${\bf R}_{\xi}^{n}$ one has the reduced problem:
\begin{equation}
\hat{u}_{tt}(t,\xi)+\vert\xi\vert^{2}\hat{u}(t,\xi) + \vert\xi\vert^{2}\hat{u}_{t}(t,\xi) = 0,\ \ \ (t,\xi)\in (0,\infty)\times {\bf R}_{\xi}^{n},
\end{equation}
\begin{equation}
\hat{u}(0,\xi)= \hat{u_{0}}(\xi),\ \ \hat{u}_{t}(0,\xi)= \hat{u_{1}}(\xi),\ \ \ x\in {\bf R}_{\xi}^{n}.
\end{equation}
Let us solve (2.1)-(2.2) directly under the condition that $0 < \vert\xi\vert \leq \delta_{0} \ll 1$. In this case we get
\[\hat{u}(t,\xi) = \frac{\hat{u_{1}}(\xi)-\sigma_{2}\hat{u_{0}}(\xi)}{\sigma_{1}-\sigma_{2}}e^{\sigma_{1}t} + \frac{\hat{u_{0}}(\xi)\sigma_{1}-\hat{u_{1}}(\xi)}{\sigma_{1}-\sigma_{2}}e^{\sigma_{2}t}\]
\begin{equation}
= \frac{e^{\sigma_{1}t}-e^{\sigma_{2}t}}{\sigma_{1}-\sigma_{2}}\hat{u}_{1}(\xi) + \frac{\sigma_{1}e^{\sigma_{2}t}-\sigma_{2}e^{\sigma_{1}t}}{\sigma_{1}-\sigma_{2}}\hat{u}_{0}(\xi),
\end{equation}
where $\sigma_{j} \in {\bf C}$ ($j = 1,2$) have forms:
\[\sigma_{1} = \frac{-\vert\xi\vert^{2}+i\vert\xi\vert\sqrt{4-\vert\xi\vert^{2}}}{2}, \quad \sigma_{2} = \frac{-\vert\xi\vert^{2}-i\vert\xi\vert\sqrt{4-\vert\xi\vert^{2}}}{2}.\]
Now let us use the idea  introduced in \cite{Ik-3} . We use the decomposition of the initial data:
\begin{equation}
\hat{u}_{j}(\xi) = A_{j}(\xi) -iB_{j}(\xi) + P_{j}\quad (j = 0,1), 
\end{equation}
where
\[A_{j}(\xi) := \int_{{\bf R}^{n}}(\cos(x\cdot\xi)-1)u_{j}(x) dx, \quad B_{j}(\xi) := \int_{{\bf R}^{n}}\sin(x\cdot\xi)u_{j}(x) dx, \quad(j = 0,1).\]
Because of (2.3) and (2.4) we get the useful identity for all $\xi$ satisfying $0 < \vert\xi\vert \leq \delta_{0}$:
\[\hat{u}(t,\xi)\ = P_{1}(\frac{e^{\sigma_{1}t}-e^{\sigma_{2}t}}{\sigma_{1}-\sigma_{2}}) + P_{0}(\frac{\sigma_{1}e^{\sigma_{2}t}-\sigma_{2}e^{\sigma_{1}t}}{\sigma_{1}-\sigma_{2}}) \]
\begin{equation}
+\,\, (A_{1}(\xi)-iB_{1}(\xi))(\frac{e^{\sigma_{1}t}-e^{\sigma_{2}t}}{\sigma_{1}-\sigma_{2}}) +  (A_{0}(\xi)-iB_{0}(\xi))(\frac{\sigma_{1}e^{\sigma_{2}t}-\sigma_{2}e^{\sigma_{1}t}}{\sigma_{1}-\sigma_{2}}).
\end{equation}
It is easy to check that
\begin{equation}
\frac{e^{\sigma_{1}t}-e^{\sigma_{2}t}}{\sigma_{1}-\sigma_{2}} = 2\frac{e^{-t\vert\xi\vert^{2}/2}\sin(\frac{t\vert\xi\vert\sqrt{4-\vert\xi\vert^{2}}}{2})}{\vert\xi\vert\sqrt{4-\vert\xi\vert^{2}}},
\end{equation}
\begin{equation}
\frac{\sigma_{1}e^{\sigma_{2}t}-\sigma_{2}e^{\sigma_{1}t}}{\sigma_{1}-\sigma_{2}} = \frac{\vert\xi\vert e^{-t\vert\xi\vert^{2}/2}\sin(\frac{t\vert\xi\vert\sqrt{4-\vert\xi\vert^{2}}}{2})}{\sqrt{4-\vert\xi\vert^{2}}} + e^{-t\vert\xi\vert^{2}/2}\cos(\frac{t\vert\xi\vert\sqrt{4-\vert\xi\vert^{2}}}{2}).
\end{equation}
If we set
\[K_{1}(t,\xi) := P_{0}\frac{\vert\xi\vert e^{-t\vert\xi\vert^{2}/2}\sin(\frac{t\vert\xi\vert\sqrt{4-\vert\xi\vert^{2}}}{2})}{\sqrt{4-\vert\xi\vert^{2}}},\]
\[K_{2}(t,\xi) := (A_{1}(\xi)-iB_{1}(\xi))(\frac{e^{\sigma_{1}t}-e^{\sigma_{2}t}}{\sigma_{1}-\sigma_{2}}), \]
\[K_{3}(t,\xi) :=  (A_{0}(\xi)-iB_{0}(\xi))(\frac{\sigma_{1}e^{\sigma_{2}t}-\sigma_{2}e^{\sigma_{1}t}}{\sigma_{1}-\sigma_{2}}),\]
then it follows from (2.5), (2.6) and (2.7) that
\[\hat{u}(t,\xi) = 2P_{1}\frac{e^{-t\vert\xi\vert^{2}/2}\sin(\frac{t\vert\xi\vert\sqrt{4-\vert\xi\vert^{2}}}{2})}{\vert\xi\vert\sqrt{4-\vert\xi\vert^{2}}} + P_{0}e^{-t\vert\xi\vert^{2}/2}\cos(\frac{t\vert\xi\vert\sqrt{4-\vert\xi\vert^{2}}}{2})\]
\begin{equation}
+\,\, K_{1}(t,\xi) +  K_{2}(t,\xi) + K_{3}(t,\xi), \quad 0 < \vert\xi\vert \leq \delta_{0}.
\end{equation}
Note that from the mean value theorem it follows that
\[2\frac{\sin(\frac{t\vert\xi\vert\sqrt{4-\vert\xi\vert^{2}}}{2})}{\vert\xi\vert\sqrt{4-\vert\xi\vert^{2}}} = \frac{2}{\sqrt{4-\vert\xi\vert^{2}}}\frac{\sin(t\vert\xi\vert)}{\vert\xi\vert} + t(\frac{\sqrt{4-\vert\xi\vert^{2}}-2}{\sqrt{4-\vert\xi\vert^{2}}})\cos(\varepsilon(t,\xi)),\]
\[\cos(\frac{t\vert\xi\vert\sqrt{4-\vert\xi\vert^{2}}}{2}) = \cos(t\vert\xi\vert) - t\vert\xi\vert(\frac{\sqrt{4-\vert\xi\vert^{2}}-2}{2})\sin(\eta(t,\xi)),\]
where
\[\varepsilon(t,\xi) := \frac{t\vert\xi\vert\sqrt{4-\vert\xi\vert^{2}}}{2}\theta + t\vert\xi\vert(1-\theta), \quad \theta \in (0,1),\]
\[\eta(t,\xi) := \frac{t\vert\xi\vert\sqrt{4-\vert\xi\vert^{2}}}{2}\theta' + t\vert\xi\vert(1-\theta'), \quad \theta' \in (0,1),\]
so that from (2.8) one has arrived at the identity:
\[\hat{u}(t,\xi) = 2P_{1}e^{-t\vert\xi\vert^{2}/2}\frac{1}{\sqrt{4-\vert\xi\vert^{2}}}\frac{\sin(t\vert\xi\vert)}{\vert\xi\vert} + P_{0}e^{-t\vert\xi\vert^{2}/2}\cos(t\vert\xi\vert)\] 
\begin{equation}
+\, P_{1}e^{-t\vert\xi\vert^{2}/2}t\frac{\sqrt{4-\vert\xi\vert^{2}}-2}{\sqrt{4-\vert\xi\vert^{2}}}\cos(\varepsilon(t,\xi)) -P_{0}e^{-t\vert\xi\vert^{2}/2}t\vert\xi\vert(\frac{\sqrt{4-\vert\xi\vert^{2}}-2}{2})\sin(\eta(t,\xi)) + \sum_{j=1}^{3}K_{j}(t,\xi).
\end{equation}
On the other hand, if one uses again the mean value theorem, it follows that 
\begin{equation}
\frac{2}{\sqrt{4-\vert\xi\vert^{2}}} = 1 + \frac{2\theta\vert\xi\vert^{2}}{(4-\theta^{2}\vert\xi\vert^{2})\sqrt{4-\theta^{2}\vert\xi\vert^{2}}}, \quad \theta \in (0,1),
\end{equation}
so that from (2.9) and (2.10) in the case when $0 < \vert\xi\vert \leq \delta_{0}$ we find that  
\[\hat{u}(t,\xi) = P_{1}e^{-t\vert\xi\vert^{2}/2}\frac{\sin(t\vert\xi\vert)}{\vert\xi\vert} + P_{0}e^{-t\vert\xi\vert^{2}/2}\cos(t\vert\xi\vert)\] 
\begin{equation}
+ P_{1}K_{4}(t,\xi) -P_{0}K_{5}(t,\xi) + K_{1}(t,\xi) +K_{2}(t,\xi) + K_{3}(t,\xi) + K_{6}(t,\xi),
\end{equation}
where
\[K_{4}(t,\xi) := e^{-t\vert\xi\vert^{2}/2}t\frac{\sqrt{4-\vert\xi\vert^{2}}-2}{\sqrt{4-\vert\xi\vert^{2}}}\cos(\varepsilon(t,\xi)),\]
\[K_{5}(t,\xi) := e^{-t\vert\xi\vert^{2}/2}t\vert\xi\vert(\frac{\sqrt{4-\vert\xi\vert^{2}}-2}{2})\sin(\eta(t,\xi)),\]
\[K_{6}(t,\xi) := P_{1}e^{-t\vert\xi\vert^{2}/2}\frac{\sin(t\vert\xi\vert)}{\vert\xi\vert}\frac{2\theta\vert\xi\vert^{2}}{(4-\theta^{2}\vert\xi\vert^{2})\sqrt{4-\theta^{2}\vert\xi\vert^{2}}}.\]

In order to prove Lemma 2.1 we have to estimate the $6$ quantities $K_{j}(t,\xi)$ ($j = 1,2,3,4,5,6$) separately. For this we prepare the following relations:
\begin{equation}
\vert\sigma_{1} - \sigma_{2}\vert = \vert\xi\vert\sqrt{4-\vert\xi\vert^{2}}, \quad (\vert\xi\vert \leq \delta_{0}),
\end{equation}
\begin{equation}
\sqrt{4-\vert\xi\vert^{2}} \geq \sqrt{4-\delta_{0}^{2}}, \quad (\vert\xi\vert \leq \delta_{0}). 
\end{equation}
\noindent
Now let us obtain several decay estimates for such $6$ quantities $K_{j}(t,\xi)$ ($j = 1,2,3,4,5,6$).\\
\noindent
{\bf (I)\,Estimate for $K_{1}(t,\xi)$.}\\ 
\[\int_{\vert\xi\vert \leq \delta_{0}}\vert K_{1}(t,\xi)\vert^{2}d\xi \leq \vert P_{0}\vert^{2}\int_{\vert\xi\vert \leq \delta_{0}}\frac{\vert\xi\vert^{2}e^{-t\vert\xi\vert^{2}}\vert\sin(\frac{t\vert\xi\vert\sqrt{4-\vert\xi\vert^{2}}}{2}) \vert^{2}}{4-\vert\xi\vert^{2}}d\xi\]
\begin{equation}
\leq \frac{\vert P_{0}\vert^{2}}{4-\delta_{0}^{2}}\int_{\vert\xi\vert \leq \delta_{0}}\vert\xi\vert^{2}e^{-t\vert\xi\vert^{2}}d\xi \leq  \frac{\vert P_{0}\vert^{2}}{4-\delta_{0}^{2}}t^{-\frac{n}{2}-1},
\end{equation}
where we have used (2.13).
\noindent
Next, we use the property that
\begin{equation}
\vert\sqrt{4-\vert\xi\vert^{2}}-2\vert = \frac{\vert\xi\vert^{2}}{2+\sqrt{4-\vert\xi\vert^{2}}} \leq \vert\xi\vert^{2}.
\end{equation}
\noindent
{\bf (II)\,Estimate for $K_{4}(t,\xi)$.}\\ 
It follows from (2.13) and (2.15) that
\[\int_{\vert\xi\vert \leq \delta_{0}}\vert K_{4}(t,\xi)\vert^{2}d\xi \leq Ct^{2}\int_{\vert\xi\vert \leq \delta_{0}}e^{-t\vert\xi\vert^{2}}\frac{\vert\sqrt{4-\vert\xi\vert^{2}}-2\vert^{2}}{(4-\vert\xi\vert^{2})}\vert\cos(\varepsilon(t,\xi)) \vert^{2}d\xi\]
\begin{equation}
\leq C\frac{t^{2}}{(4-\delta_{0}^{2})}\int_{\vert\xi\vert \leq \delta_{0}}e^{-t\vert\xi\vert^{2}}\vert\xi\vert^{4}d\xi \leq \frac{C}{(4-\delta_{0}^{2})}t^{-\frac{n}{2}}.
\end{equation}

\noindent
{\bf (III)\,Estimate for $K_{5}(t,\xi)$.}\\ 
Again it follows from (2.15) that
\[\int_{\vert\xi\vert \leq \delta_{0}}\vert K_{5}(t,\xi)\vert^{2}d\xi \leq Ct^{2}\int_{\vert\xi\vert \leq \delta_{0}}e^{-t\vert\xi\vert^{2}}\vert\xi\vert^{2}\vert \frac{\sqrt{4-\vert\xi\vert^{2}}-2}{2}\vert^{2}\vert\sin(\eta(t,\xi)) \vert^{2}d\xi\]
\begin{equation}
\leq Ct^{2}\int_{\vert\xi\vert \leq \delta_{0}}e^{-t\vert\xi\vert^{2}}\vert\xi\vert^{6}d\xi \leq Ct^{-\frac{n}{2}-1}.
\end{equation}
\noindent
{\bf (IV)\,Estimate for $K_{6}(t,\xi)$.}\\ 
\[\int_{\vert\xi\vert \leq \delta_{0}}\vert K_{6}(t,\xi)\vert^{2}d\xi \leq \vert P_{1}\vert^{2}\int_{\vert\xi\vert \leq \delta_{0}}e^{-t\vert\xi\vert^{2}}\frac{\sin^{2}(t\vert\xi\vert)}{\vert\xi\vert^{2}}\frac{4\theta^{2}\vert\xi\vert^{4}}{\vert 4-\theta^{2}\vert\xi\vert^{2} \vert^{3}} d\xi\]
\begin{equation}
\leq C\frac{\vert P_{1}\vert^{2}}{\vert 4-\delta_{0}^{2}\vert^{3}}\int_{\vert\xi\vert \leq \delta_{0}}\vert\xi\vert^{2}e^{-t\vert\xi\vert^{2}}d\xi \leq  C\vert P_{1}\vert^{2}t^{-\frac{n}{2}-1}.
\end{equation}

In order to estimate $K_{j}$ ($j = 2,3$), we prepare the following simple lemma, which plays an essential role in this research. This idea has its origin in \cite[Lemma 3.1]{Ik-3}.

\begin{lem} Let $n \geq 1$. Then it holds that for all $\xi \in {\bf R}^{n}$
\[\vert A_{j}(\xi)\vert \leq L\vert\xi\vert\Vert u_{j}\Vert_{1,1}\quad (j = 0,1),\]
\[\vert B_{j}(\xi)\vert \leq M\vert\xi\vert\Vert u_{j}\Vert_{1,1}\quad (j = 0,1),\]
where
\[L := \sup_{\theta \ne 0}\frac{\vert 1-\cos\theta\vert}{\vert\theta\vert} < +\infty, \quad M := \sup_{\theta \ne 0}\frac{\vert \sin\theta\vert}{\vert\theta\vert} < +\infty,\]
and both $A_{j}(\xi)$ and $B_{j}(\xi)$ are defined in {\rm (2.4)}.
\end{lem}
{\it Proof.} First, in the case when $\xi \ne 0$, for small $\delta > 0$ one has
\[\int_{\vert x\vert \geq \delta}\vert u_{j}(x)\vert\vert(\cos(x\cdot\xi)-1)\vert dx \leq \int_{\vert x\vert \geq \delta}\vert u_{j}(x)\vert\frac{\vert(\cos(x\cdot\xi)-1)\vert}{\vert x\cdot\xi\vert}\vert x\cdot\xi\vert dx \leq L\vert\xi\vert\int_{{\bf R}^{n}}\vert u_{j}(x)\vert\vert x\vert dx.\]
Letting $\delta \downarrow 0$ above, one has
\[\vert A_{j}(\xi)\vert \leq L\vert\xi\vert\int_{{\bf R}^{n}}\vert u_{j}(x)\vert\vert x\vert dx \quad (\xi \in {\bf R}^{n})\quad (j = 0,1).\]
Note that the above inequality holds true also in the case when $\xi = 0$. Similarly one also has
\[\vert B_{j}(\xi)\vert \leq M\vert\xi\vert\int_{{\bf R}^{n}}\vert u_{j}(x)\vert\vert x\vert dx \quad (\xi \in {\bf R}^{n}).\] 
\hfill
$\Box$
\par

\noindent
{\bf (V)\,Estimate for $K_{2}(t,\xi)$.}\,This part is crucial.\\ 
It follows from (2.6), (2.13) and Lemma 2.2 that
\[\int_{\vert\xi\vert \leq \delta_{0}}\vert K_{2}(t,\xi)\vert^{2}d\xi \leq C(L^{2}+M^{2})\Vert u_{1}\Vert_{1,1}^{2}\int_{\vert\xi\vert \leq \delta_{0}}\vert\xi\vert^{2}e^{-t\vert\xi\vert^{2}}\frac{\sin^{2}(t\vert\xi\vert\sqrt{4-\vert\xi\vert^{2}}/2)}{\vert\xi\vert^{2}(4-\vert\xi\vert^{2})} d\xi\]
\begin{equation}
\leq C\Vert u_{1}\Vert_{1,1}^{2}\int_{\vert\xi\vert \leq \delta_{0}}e^{-t\vert\xi\vert^{2}}d\xi \leq  C\Vert u_{1}\Vert_{1,1}^{2}t^{-\frac{n}{2}}.
\end{equation}

\noindent
{\bf (VI)\,Estimate for $K_{3}(t,\xi)$.}\\ 
It follows from (2.7), (2.13) and Lemma 2.2 that
\[\int_{\vert\xi\vert \leq \delta_{0}}\vert K_{3}(t,\xi)\vert^{2}d\xi \leq C(L^{2}+M^{2})\Vert u_{0}\Vert_{1,1}^{2}\int_{\vert\xi\vert \leq \delta_{0}}\vert\xi\vert^{4}e^{-t\vert\xi\vert^{2}}\frac{\sin^{2}(t\vert\xi\vert\sqrt{4-\vert\xi\vert^{2}}/2)}{(4-\vert\xi\vert^{2})} d\xi\]
\[+\,\, C(L^{2}+M^{2})\Vert u_{0}\Vert_{1,1}^{2}\int_{\vert\xi\vert \leq \delta_{0}}\vert\xi\vert^{2}e^{-t\vert\xi\vert^{2}}\cos^{2}(t\vert\xi\vert\sqrt{4-\vert\xi\vert^{2}}/2) d\xi\]
\[\leq C\Vert u_{0}\Vert_{1,1}^{2}\int_{\vert\xi\vert \leq \delta_{0}}\vert\xi\vert^{4}e^{-t\vert\xi\vert^{2}}d\xi + C\Vert u_{0}\Vert_{1,1}^{2}\int_{\vert\xi\vert \leq \delta_{0}}\vert\xi\vert^{2}e^{-t\vert\xi\vert^{2}}d\xi\]
\begin{equation}
\leq  C\Vert u_{0}\Vert_{1,1}^{2}t^{-\frac{n}{2}-2} + C\Vert u_{0}\Vert_{1,1}^{2}t^{-\frac{n}{2}-1}. 
\end{equation}

Under these preparations, let us prove Lemma 2.1.\\
{\it Proof of Lemma 2.1.} It follows from (2.11), (2.14), (2.16), (2.17), (2.18), (2.19) and (2.20) that

\[\int_{\vert\xi\vert\leq \delta_{0}}\vert \hat{u}(t,\xi)-(P_{1}e^{-t\vert\xi\vert^{2}/2}\frac{\sin(t\vert\xi\vert)}{\vert\xi\vert}+P_{0}e^{-t\vert\xi\vert^{2}/2}\cos(t\vert\xi\vert)) \vert^{2}d\xi \]
\[\leq C\vert P_{1}\vert^{2}t^{-\frac{n}{2}}+ C\vert P_{0}\vert^{2}t^{-\frac{n}{2}-1} +  C\Vert u_{1}\Vert_{1,1}^{2}t^{-\frac{n}{2}}+ C\Vert u_{0}\Vert_{1,1}^{2}t^{-\frac{n}{2}-2} + C\Vert u_{0}\Vert_{1,1}^{2}t^{-\frac{n}{2}-1}+ C\vert P_{1}\vert^{2}t^{-\frac{n}{2}-1}\]
\[\leq C\Vert u_{0}\Vert_{1,1}^{2}t^{-\frac{n}{2}-1}+ C\Vert u_{1}\Vert_{1,1}^{2}t^{-\frac{n}{2}},\]
which implies the desired estimate.
\hfill
$\Box$
\par

Based on Lemma 2.1 let us prove Theorem 1.1. For high frequency estimates we shall rely on the previous result due to \cite[Lemma 2.4]{INatsume}.\\
{\it Proof of Theorem 1.1.}
\[\int_{{\bf R}^{n}}\vert \hat{u}(t,\xi)-(P_{1}e^{-t\vert\xi\vert^{2}/2}\frac{\sin(t\vert\xi\vert)}{\vert\xi\vert}+P_{0}e^{-t\vert\xi\vert^{2}/2}\cos(t\vert\xi\vert)) \vert^{2}d\xi \]
\[= (\int_{\vert\xi\vert\leq \delta_{0}} + \int_{\vert\xi\vert\geq \delta_{0}}) \vert \hat{u}(t,\xi)-(P_{1}e^{-t\vert\xi\vert^{2}/2}\frac{\sin(t\vert\xi\vert)}{\vert\xi\vert}+P_{0}e^{-t\vert\xi\vert^{2}/2}\cos(t\vert\xi\vert)) \vert^{2}d\xi \]
\[=: I_{l}(t) + I_{h}(t).\]
To begin with, it follows from Lemma 2.1 one has
\begin{equation}
 I_{l}(t) \leq C\Vert u_{0}\Vert_{1,1}^{2}t^{-\frac{n}{2}-1}+ C\Vert u_{1}\Vert_{1,1}^{2}t^{-\frac{n}{2}}.
\end{equation}
On the other hand, it follows from \cite[Lemma 2.4 with $\theta = 1$]{INatsume} that
\[\vert\hat{u}_{t}(t,\xi)\vert^{2} + \vert\xi\vert^{2}\vert\hat{u}(t,\xi)\vert^{2} \leq Ce^{-c\rho_{\varepsilon}(\xi)t}(\vert\hat{u}_{1}(\xi)\vert^{2} + \vert\xi\vert^{2}\vert\hat{u}_{0}(\xi)\vert^{2}),\quad \xi \in {\bf R}^{n},\]
where
\[\rho_{\varepsilon}(\xi) := \varepsilon\vert\xi\vert^{2}\quad(\vert\xi\vert \leq 1), \]
\[\rho_{\varepsilon}(\xi) := \varepsilon\quad(\vert\xi\vert \geq 1),\]
where the parameter $\varepsilon > 0$ is chosen so small in the proof of \cite{INatsume}. Thus, one can estimate as follows:
\[\int_{\vert\xi\vert\geq \delta_{0}}\vert \hat{u}(t,\xi)\vert^{2}d\xi \leq C\int_{\vert\xi\vert\geq \delta_{0}}e^{-c\rho_{\varepsilon}(\xi)t}(\frac{\vert\hat{u}_{1}(\xi)\vert^{2}}{\vert\xi\vert^{2}} + \vert\hat{u}_{0}(\xi)\vert^{2})d\xi\]
\[= \int_{1 \geq \vert\xi\vert\geq \delta_{0}}e^{-c\varepsilon\vert\xi\vert^{2}t}(\frac{\vert\hat{u}_{1}(\xi)\vert^{2}}{\vert\xi\vert^{2}} + \vert\hat{u}_{0}(\xi)\vert^{2})d\xi + \int_{\vert\xi\vert\geq 1}e^{-c\varepsilon t}(\frac{\vert\hat{u}_{1}(\xi)\vert^{2}}{\vert\xi\vert^{2}} + \vert\hat{u}_{0}(\xi)\vert^{2})d\xi\]
\begin{equation}
\leq Ce^{-\alpha t}(\Vert u_{1}\Vert^{2} + \Vert u_{0}\Vert^{2}),
\end{equation}
where the constants $C > 0$ and $\alpha > 0$ depend on $\delta_{0} > 0$. On the other hand, 

\[P_{1}^{2}\int_{\vert\xi\vert\geq \delta_{0}}e^{-t\vert\xi\vert^{2}}\vert\frac{\sin(t\vert\xi\vert)}{\vert\xi\vert}\vert^{2}d\xi + P_{0}^{2}\int_{\vert\xi\vert\geq \delta_{0}}e^{-t\vert\xi\vert^{2}}\vert\cos(t\vert\xi\vert) \vert^{2}d\xi\]
\begin{equation}
\leq C(\frac{1}{\delta_{0}^{2}}+1)(\Vert u_{1}\Vert_{1}^{2} + \Vert u_{0}\Vert_{1}^{2})\int_{{\bf R}^{n}}e^{-t\vert\xi\vert^{2}}d\xi \leq C(\Vert u_{1}\Vert_{1}^{2} + \Vert u_{0}\Vert_{1}^{2})t^{-\frac{n}{2}}.
\end{equation}
Therefore, by evaluating $I_{h}(t)$ based on (2.21), (2.22) and (2.23) it is true that
\[\int_{{\bf R}^{n}}\vert \hat{u}(t,\cdot)- (P_{1}e^{-t\vert\xi\vert^{2}/2}\frac{\sin(t\vert\xi\vert)}{\vert\xi\vert}+P_{0}e^{-t\vert\xi\vert^{2}/2}\cos(t\vert\xi\vert))\vert^{2}d\xi\]
\[\leq C(\Vert u_{1}\Vert_{1}^{2} + \Vert u_{0}\Vert_{1}^{2})t^{-\frac{n}{2}} + C\Vert u_{1}\Vert_{1,1}^{2}t^{-\frac{n}{2}} + \Vert u_{0}\Vert_{1,1}^{2}t^{-\frac{n}{2}-1} + Ce^{-\alpha t}(\Vert u_{1}\Vert^{2} + \Vert u_{0}\Vert^{2}),\]
which implies the desired estimate.
\hfill
$\Box$
\par
\begin{rem}{\rm In the course of proof of Lemma 2.1 we have just encountered several singularities at $\xi = 0$ when we estimate various integrals on $\vert \xi\vert \leq \delta_{0}$ (for example, see (2.18) or (2.19) ). But, these singularities can be avoided by the same operations as in the proof of Lemma 2.2, that is, we first integrate several quantities over $[\delta, \delta_{0}]$ with sufficiently small $\delta \in (0,\delta_{0})$, and then by letting $\delta \downarrow 0$, we can have the desired estimates on $[0,\delta_{0}]$. $\xi = 0$ is the removable singularity.}
\end{rem}
\vspace{0.3cm}
\noindent{\em Acknowledgment.}
\smallskip
The work of the author was supported in part by Grant-in-Aid for Scientific Research (C)22540193 and (A)22244009 of JSPS.


\end{document}